\documentclass[12pt]{amsart}
\usepackage{pdfsync}
\usepackage{amsmath,amsthm,amsfonts,amssymb}
\usepackage{amsmath, amsthm, amssymb, amscd, amsfonts}
\usepackage[all]{xy}
\usepackage{amssymb, amsthm, amsmath}
\usepackage[english]{babel}
\usepackage{amsfonts}
\numberwithin{equation}{section}

\newcommand{\ra}{\rightarrow}

\newcommand{\Z}{\mathcal Z}
\newcommand{\ot}{\otimes}
\newcommand{\mtc}{\mathcal}

\newcommand{\Lam}{\Lambda}

\newcommand{\al}{\alpha}
\newcommand{\eps}{\epsilon}

\newcommand{\D}{\Delta}

\newcommand{\lh}{\leftharpoonup}
\newtheorem{lemma}[equation]{Lemma}
\newtheorem{thm}[equation]{Theorem}
\newtheorem{prop}[equation]{Proposition}
\newtheorem{defn}[equation]{Definition}
\newtheorem{cor}[equation]{Corollary}
\newtheorem{rem}[equation]{Remark}

\newcommand{\dw}{\downarrow}
\newcommand{\uw}{\uparrow}

\newcommand{\ch}{\chi}
\newcommand{\mtr}{\mathrm}
\newcommand{\bn}{\begin}

\newcommand{\ncm}{\newcommand}
\ncm{\np}{\newpage} \ncm{\ebl}{\end{thebibliography}}
\ncm{\bbl}{\begin{thebibliography}} \ncm{\chd}{_{ _{\ch}}}
\ncm{\ald}{_{ _{\al}}} \ncm{\cP}{\mathcal{P}} \ncm{\ei}{e_i}
\ncm{\eij}{e_{i,\;j}} \ncm{\bt}{\begin{thm}}
\ncm{\bdef}{\begin{defn}} \ncm{\edf}{\end{defn}}
\ncm{\et}{\end{thm}} \ncm{\bc}{\begin{cor}}
\ncm{\bl}{\begin{lemma}} \ncm{\el}{\end{lemma}}
\ncm{\bpf}{\begin{proof}} \ncm{\epf}{\end{proof}}
\ncm{\ec}{\end{cor}} \ncm{\er}{\end{rem}} \ncm{\br}{\begin{rem}}
\ncm{\bp}{\begin{prop}} \ncm{\ep}{\end{prop}}
\ncm{\bd}{\begin{document}} \ncm{\ed}{\end{document}}
\ncm{\beq}{\begin{equation}} \ncm{\beqn}{\begin{equation*}}
\ncm{\eeq}{\end{equation}} \ncm{\eeqn}{\end{equation*}}
\ncm{\bea}{\begin{eqnarray}} \ncm{\eea}{\end{eqnarray}}
\ncm{\beanon}{\begin{eqnarray*}} \ncm{\eeanon}{\end{eqnarray*}}
\ncm{\dimm}{\mathrm{dim}}
\ncm{\Rep}{\mtr{Rep}}\ncm{\Irr}{\mtr{Irr}}
\newcommand{\C}{\mathcal C}
\newcommand{\Di}{\mathcal D}
\newcommand{\sbst}{\subset}\thanks{This research was partially supported by CNCSIS grant PD168 no: 14/28.07.2010}
\address{Inst.\ of Math.\ ``Simion Stoilow" of the Romanian Academy
P.O. Box 1-764, RO-014700, Bucharest, Romania}\address{\textnormal{and}} \address{University of Bucharest, Faculty of Mathematics and Computer Science, 14 Academiei St., Bucharest, Romania }
\email{sebastian.burciu@imar.ro}
\date{\today}
\begin{document}
\author{Sebastian Burciu}
\title[Depth one] {Depth one extensions of semisimple algebras and Hopf subalgebras}
\maketitle
\begin{abstract}  An extension of $k$-algebras $B \subset A$ is said to have depth one if there exists a positive integer $n$ such that
$ A$ is a direct summand of  $ B^n$ in $_B\mtr{Mod} _B$.
Depth one extensions of semisimple algebras are completely characterized in terms of their centers. For extensions of semisimple Hopf algebras our results are similar to those obtained for finite group algebra extensions in \cite{BKone}.
\end{abstract}
\section{Introduction and Preliminaries}
A depth two subgroup of a finite group is
precisely a normal subgroup. This result was shown in \cite{KK}  and then extended to semisimple Hopf subalgebras in \cite{BnK}. The generalization of this result to an arbitrary extension of Hopf
algebras was then achieved in \cite{BoKu09}. A  Hopf subalgebra of a finite dimensional Hopf algebra is of depth two if and only if it is a normal Hopf subalgebra.

For an arbitrary extension of semisimple algebras a similar  to the one mentioned above also holds.  The  notion of normal Hopf subalgebras has to be replaced by a more general notion of normal subalgebras introduced in 1979 by Rieffel in \cite{Ri}.  Then the result that depth two extensions of finite dimensional semisimple algebras  coincide with normal extensions was proven in \cite{BKK09} .

Depth one extensions of semisimple algebras are in particular  depth two extensions  and therefore they are normal extensions.  However, despite of the complete understanding of depth two extensions, depth one extensions of semisimple algebras are not yet completely classified in the literature.

The only known results in this direction are those given in \cite{BKone} for group extensions.  In this paper the authors classify depth one subgroups, or equivalently depth one extensions of finite group algebras.  It is proven in \cite{BKone} that for a subgroup $H \subset G$ and an algebraically closed field $k$,  the extension $kH \subset kG$ is a depth one extension  if and only if  $H$ is a normal subgroup of $G$ and $G$ acts trivially on $\Irr(H)$.

In this paper we give a simple characterization for depth one extensions of semisimple algebras over an algebraically closed field of characteristic zero.  We show that $B \subset A$ is a depth one extension if and only if the centers of $B$ and $A$ satisfy $\Z(B)\subset \Z(A)$.

Applying this result to extensions of  semisimple Hopf algebras we obtain a similar result to the  result for group algebras from \cite{BKone} that was mentioned above.   If $B$ is a normal Hopf algebra of $A$, in Subsection \ref{action} we define an action of the fusion category $\Rep(A^*)$ on $\Rep(B)$ .  It is shown in Theorem \ref{mainhopf} that an extension $B \subset A$ of semisimple Hopf algebras is a depth one extension if and only if the fusion category $\Rep(A^*)$ acts trivially on $\Rep(B)$.

Throughout this paper $k$ is an algebraically closed field of characteristic zero.
Following  \cite{BDK} an extension of $k$-algebras $B \subset A$ is said to have left depth two (resp. right depth two) if there exists a positive integer $n$ such that
$$A \ot_B A | A^n$$ in $_B\mtr{Mod} _A$ (resp. $_A\mtr{Mod}_B$).
It is said to have depth two if it has both left depth two and right depth two. These notions of depth were introduced in \cite{KL03} and they were motivated by concepts introduced in \cite{GHJ} and \cite{KN01}. See also \cite{NV} for other motivations arising from functional analysis.

 Also an extension of $k$-algebras $B \subset A$ is said to have depth one if there exists a positive integer $n$ such that
$ A$ is a direct summand of  $ B^n$ in $_B\mtr{Mod} _B$.
It is said to have depth two if it has both left depth two and right depth two.

 We use Sweedler's notation for commultiplication  but with the sigma symbol dropped. All other Hopf algebra notations of this paper follow the standard notations  from \cite{Montg}.
\section{Depth one extensions of semisimple algebras}

Let $B \subset A$ be an extension of finite dimensional semisimple $k$-algebras.  Let $M_1, \ldots, M_s$ be  a complete set of isomorphism classes of simple $A$-modules and $e_1,\ldots ,e_s \in \Z(A)$ be their corresponding primitive central idempotents.

Similarly,  let $V_1, \ldots, V_r$  be  a complete set of isomorphism classes of simple $B$-modules and $f_1,\ldots , f_r \in \Z(B)$ be their corresponding primitive central idempotents.

Then the algebra $\mathcal{Z}(A)\cap B $ is a commutative semisimple subalgebra of both $\mathcal{Z}(A)$ and $\mathcal{Z}(B)$.
Thus there are partitions  $\mathrm{Irr}(A)=\bigsqcup_{t=1}^u\mathcal{A}_t$ and $\mathrm{Irr}(B)=\bigsqcup_{t=1}^u\mathcal{B}_t$ on the sets of irreducible representations of $A$ and $B$ such that the basis $\{m_t\}_{t=1}^u$ of primitive idempotents  of $\mathcal{Z}(A)\cap B$ can be written as
\begin{equation}
\label{mformula}
m_t=\sum_{ i \in \mathcal{A}_t}e_i=\sum_{ j \in \mathcal{B}_t}f_j.
\end{equation}
for all $1\leq t\leq u$.

The following equivalence relations on $\Irr(A)$ and $\Irr(B)$ were first introduced by Rieffel in \cite{Ri}. 
First define a symmetric relation on $\mathrm{Irr}(B)$ by $V \sim W$ if $V\uparrow^A_B$ and $W\uparrow^A_B$ have at least one  irreducible $A$-module as common constituent.
This relation $\sim$ is reflexive but not transitive in general. Its transitive closure is an equivalence relation denoted by $\approx$ or $d_B^A$.
Thus $V \approx W$ if and only if there are $V_0, \ldots, V_m \in \mathrm{Irr}(B)$ such that $V = V_0 \sim V_1 \sim \cdots \sim V_m = W$.

Also define a symmetric relation on $\mathrm{Irr}(A)$ by $M \sim N$ if $M\downarrow^A_B$ and $N\downarrow^A_B$ have a common irreducible constituent.
This relation $\sim$ is again reflexive but not transitive in general. Its transitive closure is an equivalence relation denoted by $\approx$ or $u_B^A$.
Thus $M \approx N $ if and only if there are $N_0,\ldots,N_r \in \mathrm{Irr}(A)$ such that $M=N_0\sim N_1\sim\cdots \sim N_r=N$.

The following Proposition is just a combined reformulation of Proposition 3.1 and Proposition 3.2  from \cite{BKK09}:

\begin{prop}\label{samepart} Let $B\subseteq A$ be an inclusion of semisimple finite dimensional algebras.  Then $\mathcal{A}_1,\ldots,\mathcal{A}_u$ are the equivalence classes of the equivalence relation $u_B^A$, and $\mathcal{B}_1,\ldots,\mathcal{B}_u$ are the equivalence classes of the equivalence relation $d_B^A$.  Moreover if $A$ is free as left $B$-module then with the above notations it follows that:
\begin{equation}
\label{firstformula}
[\oplus_{ M \in \mathcal{A}_t}(\dim M)\;M]\downarrow^A_B\cong\frac{\dim A}{\dim B}[\oplus_{ N \in \mathcal{B}_t}(\dim N) N]
\end{equation}
for all $1 \leq t\leq u$.
\end{prop}

In the proof of the next theorem we need the following well known result on  primitive central idempotents of semisimple algebras.  This result was also proven in \cite{BKK09}.
\begin{lemma}\label{const}Let $B \subset A$ be an extension of finte dimensional semisimple algebras.
With the above notations one has that a simple $B$-module $V_j$ is a constituent of the restriction
$M_i\downarrow^A_B$ of a simple $A$-module $M_i$ if and only if $f_j e_i\neq 0$.
\end{lemma}

Now  we can prove the following:
\bn{thm} \label{maind1}Let $B \subset A$ be an extension of  finte dimensional semisimple algebras. Then the extension $B\subset A$ is a depth one extension if and only if
$\Z(B)\subset \Z(A)$.
\end{thm}
\bn{proof}
Suppose that $B\subset A$ is a depth one extension.
Since $B$-bimodules are $B\ot B^{op}$-modules it follows that
there is some $n \geq 1$ such that $A$ is a direct summand of
$B^n$ as $B\otimes B^{op}$-modules

Let as above $e_1,\cdots , e_s$ be the primitive central idempotents of
$A$. Also let $f_1,\cdots ,f_r$ be the central primitive
idempotents of $B$.

Consider the extension of semisimple algebras $B\otimes B^{op}
\subset A\otimes A^{op}$. Then the central primitive idempotents
of $A\otimes A^{op}$ are $e_i \otimes e_j$ for all $1\leq i,j \leq
s$. Similarly the central primitive idempotents of $B\otimes
B^{op}$ are $f_i \otimes f_j$ for all $1\leq i,j \leq r$.

Note that the minimal two sided ideals of $B \otimes
B^{op}$-modules are of the form $Bf_i \otimes Bf_j$ for all $1\leq i,j \leq r$. Also $Bf_i$ is a simple $B\otimes B^{op}$-module corresponding to the minimal central idempotent $f_i \ot f_i$ and
$B=\oplus_{i=1}^sBf_j$ is a decomposition of $B$ in sum of simple $B\otimes B^{op}$-modules.

Similarly one has a decomposition $A=\oplus_{i=1}^sAe_i$ of $A$ as $A\otimes A^{op}$-modules
with $Ae_i$ simple $A\otimes A^{op}$-modules.

Since $A$ is a direct summand of $B^n$ as $B\otimes
B^{op}$-modules it follows that each simple $A\ot A^{op}$-module $Ae_i$ restricted to $B\otimes B^{op}$ is a sum of copies of the simple
$B\otimes B^{op}$- modules $Bf_j$. Recall that  the simple
$B\otimes B^{op}$- modules $Bf_j$ correspond to the central primitive idempotents $f_j\ot f_j$ of $B\ot B^{op}$. Thus the simple $B\ot
B^{op}$-modules corresponding to primitive idempotents $f_j
\otimes f_k$ with $j \neq k$ do not appear in the
decomposition of $A$ as $B\otimes B^{op}$-module.

Then Lemma \ref{const} applied to the extension $B\otimes B^{op}
\subset A\otimes A^{op}$ implies that for all $1 \leq i \leq s$ and
all $j \neq k$ one has $(f_j\ot f_k)(e_i \ot e_i)=0$. But $(f_j\ot
f_k)(e_i \ot e_i)=f_je_i\ot e_if_k$ and this implies that for all
$1\leq i \leq s$ there is a unique $j$ with $1\leq j \leq r$ such
that $e_if_j\neq 0$. Denote this index $j$ by $m(i)$. Thus for any $1 \leq i \leq s$ one has that $$e_i=e_i.1=e_i(\sum_{j=1}^rf_j)=e_if_{m(i)}.$$
On the other hand for any $j$  with $1 \leq j \leq r$ this implies that:
\begin{eqnarray*}
  f_j = f_j.1&=&f_j(\sum_{i=1}^re_i)=\sum_{\{i\;|\;f_je_i \neq 0\}}f_je_i=\sum_{\{i\;|\;m(i)=j\}}e_i.
\end{eqnarray*}
Thus $f_j$ is a central idempotent of $A$ and $\Z(B) \subset \Z(A)$.

The converse also holds.  If $\Z(B) \subset \Z(A)$ then one has $f_j
=\sum_{\{i\in S_j\}}e_i$ for all $1 \leq j \leq r$ and for some pairwise disjoint subsets $S_j \subset \{1,\cdots ,\;s\}$. Thus for all $1 \leq i \leq s$ the simple constituents of  the restriction of
$Ae_i$ as a $B\ot B^{op}$-module are those of the type $Bf_j$
 where $j$ is chosen such that $i \in S_j$. This implies
there is $n$ such that $A$ is a direct summand of $B^n$ as
$B\otimes B^{op}$-modules.
\end{proof}
\section{Depth one extensions of semisimple Hopf algebras}\label{d1hopf}
\subsection{Preliminaries on semisimple Hopf algebras}
Let $A$ be a finite dimensional semisimple Hopf algebra over $k$. Then $A$ is also cosemisimple \cite{Lard}. Denote by $\mtr{Irr}(A)$ the set of irreducible characters of $A$ and by $C(A)$ the character ring  of $A$. Then $C(A)$ is a semisimple subalgebra of $A^*$ \cite{Z} and $C(A)=\mtr{Cocom}(A^*)$,
the space of cocommutative elements of $A^*$. By duality,
the character ring of $A^*$ is a semisimple subalgebra of $A$ and under this
identification it follows that $C(A^*)=\mtr{Cocom}(A)$. 

If $B$ is a Hopf subalgebra of $A$ recall that $B$ is called a normal Hopf subalgebra if $B$ is closed both under left and right adjoint action of $A$ on itself. That is $a_1BS(a_2)\subset B$ and $S(a_1)Ba_2\subset B$ for all $a \in A$. If $A$ is a semisimple  Hopf algebra then $B$ is a normal Hopf subalgebra of $A$ if and only if $a_1BS(a_2)\subset B$ for all $a \in A$.
\subsection{Formulas for the central primitive idempotents}
Let $B\subset A$ be an inclusion of semisimple Hopf algebras over
an algebraically closed field $k$.  As before let $e_1,\cdots
e_r$ be the central idempotents of $A$ and let $\ch_1 ,\cdots ,
\ch_r$ be the characters associated to them.  The following formula for  primitive central idempotents is well known:
\bn{equation}\label{a}
e_i=\frac{\ch_i(1)}{\dim A}(\Lam \lh S(\ch_i))=\frac{\ch_i(1)}{\dim A}\ch_i(S(\Lam_1))\Lam_2
\end{equation}
where $\Lam$ is the idempotent integral of $A$.

Similarly let $f_1, \cdots ,f_s$ be the central idempotents of $B$
and $\al_1, \cdots ,\al_s$ be the corresponding  irreducible characters of $B$. A
formula similar to the one above for the idempotents $f_j$ of $B$ gives that:
\bn{equation}\label{b}
f_j=\frac{\al_j(1)}{\dim B}(\Gamma \lh S(\al_j))=\frac{\al_j(1)}{\dim B}\al_j(S(\Gamma_1)\Gamma_2
\end{equation}
where $\Gamma$ is the idempotent integral of $B$.

\subsection{ Simple subcoalgebras associated to a simple comodule and a formula for the integral}
 Let $W$ be a simple left $A^*$-module. Then $W$ is a simple right $A$-comodule and  one can associate to it a simple subcoalgebra of $A$ denoted by $C_{_W}$ \cite{Lar}. If $q=|W|$ then $|C_{_W}|=q^2$ and it is a matrix coalgebra. It has a basis $\{x_{ij}\}_{1\leq i,j \leq q}$  such that $\D(x_{ij})=\sum_{l=1}^qx_{il}\ot x_{lj}$ for all $1\leq i,j \leq q$. Moreover, the character of $W$ as left $A^*$-module is $d \in C(A^*)\subset A$ and is given by $d=\sum_{i=1}^qx_{ii}$. Then $\eps(d)=q$ and the simple subcoalgebra $C_{ _W}$ is also denoted by $C_d$.

It follows from Proposition 4.1  of  \cite{Lar} that the idempotent integral $\Lam$ of $A$ can be written as:
\beq\label{intform}
\Lam=\frac{1}{\dim A}\sum_{d \in \mtr{Irr}(A^*)}\eps(d)d.
\eeq

\subsection{Conjugate modules and the action of $\Rep(A^*)$ on $\Rep(B)$}\label{action}

The following results are taken from \cite{cos}. Let $B$ be a normal Hopf subalgebra of $A$ and $M$ be a $B$-module.  If $W$ is an $A^*$-module then $W\ot M$ becomes a $B$-module with
\bn{equation}\label{def}
 b(w\ot m)=w_0\ot(S(w_1)bw_2)m,
\end{equation}
for all $b \in B$, $w \in W$ and $m \in M$.

It can be checked that if $W \cong W'$ as $A^*$-modules then $W\ot M \cong W' \ot M$. Thus for any irreducible character $d \in \mtr{Irr}(A^*)$ associated to a simple $A$-comodule $W$ one can define the $B$-module $^dM \cong W \ot M$. In analogy with the group  algebra situation the module $W\ot M$ is called the conjugate module of $M$ by $W$.

It can be checked that this defines an action (not necessarily monoidal in general) of $\Rep(A^*)$ on $\Rep(B)$. For the definition of  the action of a monoidal category on an abelian category one can consult \cite{d2}.

We say that the action of  $\Rep(A^*)$ on $\Rep(B)$ is trivial if $W\ot M \cong (\dim W)M$ for any $A$-comodule $W$ and any $B$-module $M$.

\bn{prop} \label{cojre}(Proposition 4.19, \cite{cos}.) Let $B$ be a normal Hopf subalgebra of a semisimple Hopf algebra $A$ and $M$ be an irreducible $B$-module. Then $M \uw_{ _B}^A\dw_{ _B}^A$ and $\oplus_{d \in \mtr{Irr}(A^*)}\;^dM$ have the same irreducible constituents.
\end{prop}

It follows from the above Proposition that   $\Rep(A^*)$ acts trivially on $\Rep(B)$ if and only if
\beq\label{trivact}
M \uw_{ _B}^A\dw_{ _B}^A=\frac{\dim A}{\dim B}M
\eeq
for any $B$-module $A$.
\subsection{Depth one extensions of semisimple Hopf algebras}

In this subsection we give a characterization of depth one extensions of semisimple Hopf algebras.  It will be shown the following:

\bt\label{mainhopf}
 Suppose that $B$ is a Hopf subalgebra of $A$.  Then $B \subset A$ is a depth one extension of semisimple Hopf algebras if and only if $B$ is a normal Hopf subalgebra of $A$ and $Rep(A^*)$ acts trivially on $Rep(B)$.
\et
\bpf Suppose that $B \subset A$ is a depth one extension of semisimple Hopf algebras.  Then $B \subset A$ is also a depth  two extension and by Theorem 2.10 from \cite{BoKu09} it follows that $B$ is a normal Hopf subalgebra of $A$.  Also by  Theorem \ref{maind1} above it follows that $\Z(B)\subset \Z(A)$ and therefore
\beq
f_j=\sum_{\{i\;|\;f_je_i \neq 0\}}e_i.
\eeq
As before, denote by $S_j$ the set of all indices $i$  with $f_je_i \neq 0$. Applying formula \ref{a}  for each $e_i$ it follows that
\beq
f_j=\frac{1}{\dim A}\sum_{d \in \mtr{Irr}(A^*)}<\sum_{i\in S_j}\ch_i(1)\ch_i,\;S(x_{uv}^d)>x_{vu}^d.
\eeq

On the other hand from formula \ref{b} it follows that
\beq
f_j=\frac{1}{\dim B}\sum_{d \in \mtr{Irr}(B^*)}<\al_j(1)\al_j,\;S(x_{uv}^d)>x_{vu}^d.
\eeq

Comparing the last two formulae it follows that
\beq
<\frac{1}{\dim A}\sum_{i\in S_j}\ch_i(1)\ch_i,\;S(x_{uv}^d)>=\frac{1}{\dim B}<\al_j(1)\al_j,\;S(x_{uv}^d)>
\eeq

for all $d \in\mtr{Irr}(B^*)$ and
\beq
<\frac{1}{\dim A}\sum_{i\in S_j}\ch_i(1)\ch_i,\;S(x_{uv}^d)>=0
\eeq
if $d \notin \mtr{Irr}(B^*)$.

This shows that
\beq\label{x}
(\sum_{i\in S_j}\ch_i(1)\ch_i)\dw^A_B=\frac{\dim A\;\al_j(1)}{\dim B}\al_j
\eeq
for all $1 \leq j \leq r$.
Note that by Frobenius reciprocity  the above relation \ref{x} implies that
\beq\label{y}
\al_j\uw^A_B=\frac{\al_j(1)\dim B}{\dim A}(\sum_{i\in S_j}\ch_i(1)\ch_i).
\eeq
This shows that  $\al_j\uw^A_B\dw^A_B=\frac{\dim A}{\dim B}\al_j$ for all irreducible characters $\al_j$. Then  the remark given after Proposition \ref{cojre} implies that $\Rep(A^*)$ acts trivially on $\Rep(B)$.

The converse of the theorem is immediate. 
Suppose that $B$ is a normal  Hopf subalgebra of $A$ and  $\Rep(A^*)$ acts trivially on $\Rep(B)$. Then Proposition 5.8 of \cite{cos} implies that $M\dw^A_B$ is an homogenous $B$-module for any simple $A$-module $M$. Applying Proposition  \ref{samepart} it follows that  for any $1 \leq t \leq u$ each equivalence class $\mtc{B}_t$ has just one element.  Then formula \ref{mformula} implies that $\Z(B)\subset \Z(A)$ and therefore $B\subset A$ is a depth one extension by Theorem \ref{maind1}.
\epf
 For any subcoalgebra $C$ of $A$ let $C\otimes B$ be the $B\ot B^{op}$-module with the structure given by $$b(c \ot x)b'=c_1\ot (Sc_2bc_3)xb'$$
for all $c \in C$ and $b,x,b'\in B$. It is easy to check that $C\ot B$ is a $B\ot B^{op}$-module with the above structure.
\bc\label{final}
Let $B \subset A$ be a  Hopf subalgebra of a semisimple Hopf algebra  $A$.  Then the following are equivalent:

1) $B$ is a depth one Hopf subalgebra of $A$.

2) $A \cong B^{[A:B]}$ in $_B\mtr{Mod}_B$.

3) $B$ is a normal Hopf subalgebra of $A$ and one has $C\ot B \cong B^{\dim\; C}$ in $_B\mtr{Mod}_B$, for any simple subcoalgebra $C$ of $A$ .

4) $B$ is a normal Hopf subalgebra of  $A$ and $\Rep(A^*)$ acts trivially on $\Rep(B)$.
\end{cor}

\bpf
$1)\Rightarrow 2)$ Suppose that $B$ is a depth one subalgebra of $A$. Then $\Z(B)\subset \Z(A)$  and with the notations from the proof of Theorem \ref{maind1} write $$f_j=\sum_{i\in S_j}e_i$$ for all $1\leq j \leq r$.
Thus $Ae_i$ regarded as $B$-module is a homogenous module with the isomorphism type of the simple $B$-module corresponding to the unique primitive central idempotent $f_j$ of $B$ with $i \in S_j$. Since $A$ is free as left $B$-module it follows that $\oplus_{i\in S_j}Ae_i=(Bf_j)^{[A:B]}$ as left $B$-modules. Therefore
$$
\sum_{i \in S_j}\dim Ae_i=[A:B](\dim Bf_j)
$$
On the other hand from the proof of Theorem \ref{maind1} it follows that for any $i$ with $1 \leq i \leq s$ there is a unique $j$ such that
$$(Ae_i)\dw^{A\otimes A^{op}}_{B\ot B^{op}}=(Bf_j)^{\frac{\dim Ae_i}{\dim Bf_j}}.$$
Thus \bn{eqnarray*}
A & = & \oplus_{i=1}^r(Ae_i)\\ &=&
\oplus_{j=1}^s\oplus_{i\in S_j}(Ae_i)
\\ &=& \oplus_{j=1}^s(Bf_j)^{[A:B]}\\ &=&B^{[A:B]}.
\end{eqnarray*}
as $B$-bimodules.

$2)\Rightarrow 1)$ clear.

$1)\Rightarrow 3)$   Suppose that $B$ is a depth one Hopf subalgebra of $A$.  Then as above $B$ is a normal Hopf subalgebra of $A$. On the other hand  $B\ot B^{op}$ is also  a depth one Hopf subalgebra of $A\otimes A^{op}$. Indeed, $\Z(B\ot B^{op})=\Z(B)\ot \Z(B)\subset \Z(A)\otimes \Z(A)=\Z(A\otimes A^{op})$ since $\Z(B)\subset \Z(A)$  by Theorem \ref{maind1} . Apply Theorem \ref{mainhopf} for the extension  $B\ot B^{op} \subset A\ot A^{op}$.  Take $C\ot k$ as a simple subcoalgebra of $A\ot A^{op}$ and $W:=B$ as a $B\ot B^{op}$-module. We claim that the conjugate module of $W$ by $C\ot k$ is the $B$-bimodule $C\otimes B$ defined above.  Indeed as vector spaces the conjugate module is $C\ot k\ot B$ and the strucure is given by $$(b\ot b')(c
\ot 1 \ot x)=(c_1\ot 1)\ot (Sc_2\ot 1)(b\ot b')(c_3\ot 1)x=(c_1\ot 1)\ot(Sc_2bc_3)xb'$$
for all $c \in C$ and $b,x,b'\in B$.  Then Theorem \ref{mainhopf} implies that $C\ot B \cong B^{\dim C}$ as $B\ot B^{op}$-modules.

$3)\Rightarrow 1)$ Recall the coset decomposition $$A=\oplus_{C \in \mtr{Irr}(A^*)/\sim }CB$$ as $B\ot B^{op}$-modules from  Corrolary 2.5  of \cite{cos}. Note that each $CB$ is a $B\ot B^{op}$-submodule of $A$ since $b(cx)b'=c_1(Sc_2bc_3)xb'\in CB$ for all $c \in C$ and $b,x,b'\in B$.

It is easy to verify that the multiplication map $C\ot B \ra CB$ is a morphism of $B$-bimodules. Since $B\ot B^{op}$ is semisimple it follows that each $B$-bimodule $CB$ is a direct summand of the bimodule $B^{\dim C}$. This implies that $A=\oplus_{C \in \mtr{Irr}(A^*)/\sim }CB$ is a direct summand of $B^n$ as $B\ot B^{op}$-modules for some $n \geq 0$.

$1)\Leftrightarrow 4)$ is Theorem \ref{mainhopf} from above.
\epf

\bibliographystyle{amsplain}
\bibliography{fusion}
\end{document}